\documentclass[12pt]{article}
\usepackage{amsmath,amssymb,amsthm,amsfonts,graphicx,color,enumerate,float,epsfig}
\begin{document}

\newcommand{\kf}[1]{\marginpar{\tiny #1 --kf}}
\newtheorem{thm}{Theorem}[section]
\newtheorem{lemma}[thm]{Lemma}
\newtheorem{cor}[thm]{Corollary}
\newtheorem{prop}[thm]{Proposition}
\newtheorem{question}[thm]{Question}
\newtheorem*{mainthm0}{Theorem \ref{main}}
\newtheorem*{mainthm1}{Corollary \ref{manifold}}

\newtheorem{remark}[thm]{Remark}
\newtheorem{important remark}[thm]{Important Remark}
\newtheorem{definition}[thm]{Definition}
\newtheorem{example}[thm]{Example}
\newtheorem{fact}[thm]{Fact}
\newtheorem{convention}[thm]{Convention}

\def\diam{\operatorname{diam}}
\def\vl{\operatorname{vl}}
\def\svl{\operatorname{svl}}
\def\cal{\mathcal}
\def\R{{\mathbb R}}
\def\Z{{\mathbb Z}}
\def\D{\Delta}

\def\Mod{{\rm Mod}}
\def\ep{{\varepsilon}}
\def\Isom{{\rm Isom}}
\def\C{{\mathcal C}}
\def\T{{\mathcal T}}
\def\TT{\overline{\mathcal T}}
\def\dTT{\partial \overline{\mathcal T}}
\def\S{{\mathcal S}}
\def\SS{{\overline{\mathcal S}}}
\def\M{{\mathcal M}}
\def\L{{\mathcal L}}
\def\D{\Delta}
\def\s{\sigma}
\def\t{\tau}
\newcommand{\cN}{{\cal N}}

\def\red{\textcolor{red}}

\title{Asymptotically isometric  
metrics on relatively 
hyperbolic groups and marked length spectrum}

\author{Koji Fujiwara\thanks{The author is supported in part by
Grant-in-Aid for Scientific Research (No. 23244005)}
} 


\maketitle

\begin{abstract}
We prove asymptotically isometric, coarsely geodesic metrics
on a toral relatively hyperbolic  group
are coarsely equal. The theorem applies to all lattices in 
$SO(n,1)$. This partly verifies a conjecture by Margulis.
In the case of hyperbolic groups/spaces, our result generalizes
a theorem by Furman and a theorem by Krat. 

We discuss  an application to the isospectral
problem for the length spectrum of  Riemannian manifolds.
The positive answer to this problem has been known for several cases.
All of them have hyperbolic fundamental groups.
We do not solve the isospectral problem
in the original sense, but prove the universal covers
are $(1,C)$-quasi-isometric if the fundamental 
group is a toral relatively
hyperbolic group.

\end{abstract}

\section{Introduction}
\subsection{Asymptotically isometric metrics}
Suppose a group $G$ acts on a space $X$ with two metrics
$d_1$ and $d_2$ that are $G$-invariant. 
We say $(X,d_1)$ and $(X,d_2)$ are 
\begin{enumerate}
\item
{\it coarsely equal} if there exists $C$ such that 
for all $x,y \in X$,
$$|d_1(x,y) -d_2(x,y) | \le C.$$
\item
{\it asymptotically isometric} (or jsut {\it asymptotic}) if 
$d_1(x,y) \to \infty$ if and only if $d_2(x,y) \to \infty$, and 
$$\frac{d_1(x,y)}{d_2(x,y)} \to 1 \,\,\,  {\rm as} \,\,\, d_2(x,y) \to \infty.$$
\item
{\it weakly asymptotically isometric} (or {\it weakly asymptotic}) if
$(X,d_1)$ and $(X,d_2)$ are quasi-isometric, $g \in G$
is hyperbolic for $d_1$ if and only if it is hyperbolic for $d_2$, and 
that 
for every $g \in G$ that is hyperbolic,
we have 
$$\lim _{n \to \infty} \frac{d_1(x, g^n(x))}{d_2(x,g^n(x))}=1.$$
This property does not depend on the choice of $x$.
\item 
$d_1$ and $d_2$  have the {\it same marked length spectrum}
with respect to the $G$-actions 
if 
 $|g|_1=|g|_2$ for all $g \in G$,
 where 
$|g|_i$ is the {\it translation length} of $g$ for $d_i$ defined by
$$|g|_i= \lim_{n\to \infty} \frac{d_i(x,g^n(x))}{n}.$$
$g$ is {\it hyperbolic} on $(X,d_i)$ if and only if $|g|_i>0$.
\end{enumerate}

Clearly, (1) $\Rightarrow$ (2) $\Rightarrow$ (3)
$\Leftrightarrow$ (4). We discuss
the other implications in this paper. 
Some remarks are in order.
(1) is same as the identity map is a $(1,C)$-quasi-isometry.
This is stronger than $d_1$ and $d_2$ are $(1,C)$-quasi-isometric.
Even if $d_1$ and $d_2$ are isometric, (1) may not hold. 
(2) implies that the identity map is a quasi-isometry.
In view of that, if we concern the implication from (3) to 
(1) or (2), then we should look at $G$-equivariant maps 
from $(X,d_1)$ to $(X,d_2)$, not 
the identity map.

%

We write $A \sim _C B$ if $|A-B|\le C$.
A metric space $(X,d)$ is {\it $C$-coarsely geodesic} if
for any $x,y \in X$, there is a path $\gamma$ from 
$x$ to $y$ such that for all $t, s$, we have 
$|t-s| \sim _C |\gamma(t)-\gamma(s)|$.
$\gamma$ does not have to be continuous.
If $C=0$, $X$ is geodesic.
We may suppress the constant $C$ and just say {\it coarsely geodesic}.

We write $d_1 \sim _C d_2$ if 
for all $x,y$ we have $d_1(x,y)  \sim _C d_2(x,y)$.
We may simply write $d_1 \sim d_2$.
This is nothing but they are coarsely equal.



Recall that a map $f : X \to Y$ between two metric spaces $(X, d_X)$ and $(Y, d_Y )$ is
called a $(L,C)$-{\it quasi-isometry} if $Ld_X(a, b) - C \le d_Y (f(a), f(b)) \le d_X(a, b)/L + C$, for
all $a, b \in X$, and every $y \in Y$ is at distance at most $C$ from some element of $f(X)$. If such $f$ exists for some $L \ge 1, C\ge 0$, then $X$ and $Y$ are quasi-isometric.

Burago \cite{burago} proved for $G={\Bbb Z}^n$
and $X=\R^n$, (2) $\Rightarrow$ (1) for  $G$-invariant Riemannian metrics.
His argument applies to a pair of coarsely geodesic 
metrics  on ${\Bbb Z}^n$ (Corollary \ref{burago2}).

Krat \cite{krat} proved an analogous result
when $X$ is $\delta$-hyperbolic, in particular, which implies 
(2) $\Rightarrow$ (1) for two left invariant metrics
on a hyperbolic group $G$ that are quasi-isometric to a word metric.
Furman \cite{furman} proved (3) $\Rightarrow$ (1) in the same setting.
His argument is different from hers.
We will modify her argument and prove (3) $\Rightarrow$ (1)
for toral relatively hyperbolic groups, which are more general 
than hyperbolic groups.

Abels and Margulis \cite{abels} proved (2) $\Rightarrow$ (1)
for ``word metrics" on reductive Lie groups.
It is asked by Margulis in \cite{margulis} if
(2) $\Rightarrow$ (1) holds in general on a finitely generated group $G$.
Breuillard \cite{breu} answered this question in the negative.
He found a counter example, which are two word-metrics 
on $H_3(\Bbb Z) \times \Bbb Z$, where $H_3(\Bbb Z)$ is the three dimensional
discrete Heisenberg group.
In fact, the two metrics are  not coarsely equal on some cyclic (undistorted) subgroup
in his example. 
Also, those two metrics are 
not even $(1,C)$-quasi-isometric to each other for any $C$, \cite{breu2}.
It was known that (2) $\Rightarrow$ (1)  on $H_3(\Bbb Z)$, \cite{krat}.

\subsection{Main results}

An isometric action on a metric space $(X,d)$ by a group $G$
is {\it cobounded} if there exists a bounded set $B$ in $X$ such that 
$G.B=X$, and {\it proper}  if for any $x \in X$ and $R>0$
there exist at most finitely many $g \in G$ with 
$d(x, g.x) \le R$.

The following is the main result. 
\begin{mainthm0}
Assume $(G, \mathcal H)$ is a  relatively hyperbolic group such that 
 for each $H_i \in \mathcal H$, $H_i$ contains $\Z^{n_i}$ as a finite index subgroup.
Assume $G$ acts on $X$ properly and co-boundedly by isometries
for geodesic metrics $d_1, d_2$ (or more generally, $d_2$
is a coarsely geodesic metric). 
If they are weakly asymptotically isometric, then 
they are coarsely equal.
\end{mainthm0}


The following are examples of toral relatively hyperbolic groups
(see \cite{hruska} and Theorem 1.2.1 therein):
\begin{itemize}
\item
all lattices in $SO(n,1)$ (uniform ones are hyperbolic).
\item
CAT(0) groups with isolated flats (\cite{hruska}).
In particular, the fundamental group of 
a closed, irreducible 
 $3$-manifold such that each piece of its JSJ-decomposition
is atoroidal (namely, hyperbolic).

\item
Limit groups in the sense of Sela.

\end{itemize}
It seems it is an open question 
if the conclusion of the theorem holds for 
 non-uniform lattices
in the Lie group $SU(n,1)$.  See the discussion 
in Section \ref{complex}.


A merit to show (3) $\Rightarrow$ (1) is  it has 
 an application to the marked length spectrum problem
since (3) $\Leftrightarrow$ (4).

\begin{mainthm1}

Let $(M_1,d_1), (M_2,d_2)$ be closed Riemannian manifolds
with the isomorphic fundamental group $G$ that 
is toral relatively hyperbolic.
Assume they have the same marked length spectrum.
Then there is a $G$-equivariant $(1,C)$- quasi-isometry map
$f:\tilde M_1 \to \tilde M_2$.
\end{mainthm1}

It is easy to see that $M_1$ and $M_2$ have the 
same marked length spectrum if the conclusion in the corollary holds. 
Therefore we rephrase the marked length spectrum problem as follows.
If there is a $G$-equivariant $(1,C)$- quasi-isometry map
from $\tilde M_1$ to $\tilde M_2$, then is $M_1$ isometric to $M_2$ ?
Notice that if $C=0$, then $M_1$ and $M_2$ are isometric.
 
The iso-spectral problems for the marked length spectrum 
has been solved for several families of Riemannian manifolds, but 
in all of those cases, the fundamental 
group is hyperbolic (see Section \ref{marked}).
The novel part of our result is that 
we put  the marked length spectrum 
problem into context for a broader class
of groups.

We close the introduction with a discussion on a question by Gromov.
In \cite[2C$_2$(c)]{gromov.asym} he asks
if the Hausdorff distance of two manifolds $X_1$ and $X_2$ 
is finite if they are acted by $G$ properly and co-compactly
and that ``$AL \, {\rm Dist}(X_1,X_2)=0$". Here, AL is for asymptotically 
Lipschitz and $AL \, {\rm Dist}(X_1,X_2)=0$
if for any $a>0$ there is a $(1+a, C_a)$-quasi-isometry 
$f_a$ between $X_1$ and $X_2$ for some $C_a$.
Our property (2) implies $AL \, {\rm Dist}(X_1,X_2)=0$
since we can take the identity map as $f_a$ for any $a>0$
with sufficiently large $C_a$.
(1) implies that the Hausdorff distance 
between $X_1$ and $X_2$ is at most $C$ via the identity map.
Our result affirmatively answers a version of the question by Gromov with both the
assumption and the conclusion $G$-equivariant  when $G$ is toral 
relatively hyperbolic.

\section{The case of hyperbolic spaces}

We first prove the following.
 Krat \cite{krat} proved the result under the assumption
that $d_1, d_2$ are asymptotic, but 
our assumption is weaker. 

\begin{thm}\label{krat}
Let $d_1, d_2$ be coarsely 
geodesic metrics on $X$ on which $G$ acts
by isometries, co-boundedly with respect to both $d_1$ and $d_2$.
Suppose $(X,d_1)$ is $\delta$-hyperbolic.
Assume that $d_1$ and $d_2$ are weakly 
asymptotically isometric.
Then, $d_1 \sim d_2$.
\end{thm}


\proof
The outline of our argument is same as the one by Krat. 
Define $\D(x,y)=d_1(x,y)-d_2(x,y)$.
To argue by contradiction, assume $\D$ is not bounded. We will prove  that 
then $\liminf_{n} \D(x,g^n(x))/d_2(x,g^n(x)) \not= 0$ as $n \to \infty$
for some $g \in G$.
This contradicts to the assumption.

Here are two elementary lemmas.
The first one is straightforward
from the triangle inequality for $d_1$ and $d_2$.
\begin{lemma}\label{triangle}
Let $x,y,z$ be points in $X$. Then,

$|\D(x,y) - \D(x,z) |  \le d_1(y,z) +d_2(y,z)$.
\end{lemma}

\begin{lemma}\label{basic.additive}
Let $\gamma_1$ be a $d_1$-geodesic and $\gamma_2$ a $d_2$-geodesic
from $x$ to $y$. Let $z \in \gamma_1$ be a point such that 
there exists $z' \in \gamma_2$ with $d_2(z,z')\le C$.
Then,
$\D(x,z)+\D(z,y) \sim _{2C} \D(x,y)$.
Moreover, the conclusion holds if $\gamma_1$ and $\gamma_2$
are $L$-coarse geodesics (with a constant that is larger than $2C$
depending on $L$).

\end{lemma}
We say that $\D$ is {\it almost additive at $z$}.
\proof
Suppose $\gamma_1, \gamma_2$ are geodesics.
$\D(x,z)=d_1(x,z)-d_2(x,z) \sim_{C} d_1(x,z) -d_2(x,z')$,
and similarly $\D(z,y) \sim_{C} d_1(z,y)-d_2(z',y)$.
Now, $\D(x,z)+\D(z,y) \sim_{2C}
d_1(x,z) -d_2(x,z') +(d_1(z,y)-d_2(z',y))=d_1(x,y)-d_2(x,y)=\D(x,y).$

A similar argument applies when $\gamma_i$ are coarse geodesics
and we omit details. 
\qed

We go back to the proof of the theorem. In the following 
we assume that $d_1$ and $d_2$ are geodesic metrics.
We can easily modify each argument when they are only coarsely
geodesic with extra constants, but we leave it to the readers.

Since $d_1$ and $d_2$ are asymptotic, therefore quasi-isometric to each other, 
any $d_1$-geodesic is a $d_2$-quasi-geodesic with a controlled 
quasi-isometric constants, vice-versa.
There exists a constant  $C$ (by  the Morse lemma, \cite[Theorem 1.7]{bridson})
such that a $d_1$-geodesic and a $d_2$-geodesic with the 
same endpoints are in the $C$-neighborhood (for both $d_1$ and $d_2$) of each other.
Therefore by Lemma \ref{basic.additive},
$\Delta$ is almost additive, for a uniform constant, on any $d_1$-geodesic at any point.


By assumption, a metric ball of radius, say, $D$ covers $X$ by the $G$-action. 
Fix a base point $x$ and we write   $g(x)$ as $g$.
Write $\Delta(1,g)$ as $\D(g)$. Notice $\D(g)=\D(g^{-1})$.
By assumption there is $g$ with $\D(g) >> \delta, C, D$.
Take a $d_1$-geodesic $\gamma$ from $1$ to $g$, then 
there is $h\in G$ with $\D(h)$ is approximately $\D(g)/2$.
This is possible since $\D$ is almost continuous on a geodesic. 
Set $k=h^{-1}g$. Then 
$\D(k)$ is approximately $\D(g)/2$ by Lemma \ref{basic.additive}
and Lemma \ref{triangle}.

Let $[p,q]$ denote a $d_1$-geodesic from $p$ to $q$.
For $g \in G$, let $g^*.x$ denote the  piecewise 
geodesic $\cup_{n \ge 0}[g^n.x, g^{n+1}.x]$ starting at $x$.
\\
{\bf Claim}. At least one of the following piecewise
geodesic is a
$(2,10\delta+D)$-quasi-geodesic on $X$:
$k^*.x$ or  $h^*.x$ or $(hk)^*.x$.


This is a standard fact. See Lemma 8.1.A in \cite{gromov}.
Since the points $1,h,g=hk$ is almost on a geodesic,
 $d_1(x,g.x) \sim d_1(x,h.x) + d_1(h.x,g.x)$ and 
both $d_1(x,h.x)$ and $d_1(h.x, g.x)$ are large,
since $\Delta(g)$ is large.
This is enough to apply Lemma 8.1.

This claim is the new ingredient than the argument by Krat, who stated
that a piecewise geodesic of a similar property, not necessarily periodic,
exists (Lemma 1.12).
The periodicity is 
crucial  under our assumption. 

Let $f^*.x$ be one of the paths we obtain from the claim.
On the path, $\D(x,f^n.x)$ grows roughly linearly on $n$
since on each geodesic piece it increases at least by $\D(h)-2C$ or $\D(g)-2C$,
depending on $f^*.x$, by Lemma \ref{basic.additive}. 
We used that $\D(h) \sim \D(k) >> \delta, C, D$.
It follows $\liminf_n \D(x,f^n(x))/d_1(x,f^n(x)) >0$.
\qed

%
%

\section{Toral relatively hyperbolic groups}

In this section we will generalize Theorem \ref{krat}
to relatively hyperbolic groups. 

\subsection{Asymptotically tree graded spaces and relatively hyperbolic groups}
We review a few key notions and results from \cite{DS}.
Let $X$ be a complete geodesic metric space and let $\mathcal P$ be a collection 
of closed geodesic subsets (called {\it pieces}) with the following properties
(\cite{DS}):
\\
(T1) every two different pieces have at most one common point.
\\
(T2) Every simple geodesic triangle (a geodesic triangle that is 
a simple loop) in $X$ is contained in one piece.

Let $X$ be a metric space and $\mathcal A$ a collection of subsets in $X$.
$X$ is {\it asymptotically tree-graded with respect to $\mathcal A$}
if every asymptotic cone, $Con^{\omega}(X)$, of $X$ is tree-graded with respect
to a certain collection of subsets $\mathcal A_{\omega}$, which are defined from the 
collection $\mathcal A$ (see \cite[Definition 3.19]{DS} for the precise definition).
Here, $\omega$ is an ultra filter. See the definition of the asymptotic 
cone in \cite[Definition 3.8]{DS}.
In this paper, we do not use the definitions of (asymptotically) tree-graded spaces,
but only quote geometric properties of those spaces from \cite{DS}
and \cite{sisto}. We will state them later. 

There is more than one way to define relatively hyperbolic groups.
The following definition is one of the main theorems in \cite{DS}.
We say a finitely generated group $G$ is {\it relatively hyperbolic} with respect
to a collection of subgroups $\mathcal H=\{H_1, \cdots, H_m\}$ if 
$G$ is asymptotically tree-graded with respect to the subgroups
$\mathcal H$, namely, 
the Cayley graph of $G$, $\Gamma(G,S)$, with respect to some (and any)
finite set, $S$, of generators is asymptotically tree-graded
with respect to the collection of left cosets $\{gH_i | g \in G, i=1, \cdots, m\}$.
The subgroups $H_i$ are called {\it peripheral} subgroups.
If they are all finitely generated virtually abelian groups, $G$ is 
called a {\it toral} relatively hyperbolic group.
We do not 
assume that a toral relatively hyperbolic group is torsion-free.

Farb \cite{farb} defined that $G$ is {\it weakly relatively hyperbolic} with respect to 
$\mathcal H$ if the Cayley graph $\Gamma(G,S \cup \mathcal H)$ is
hyperbolic, where $S \cup \mathcal H$
means the union of the elements in $S$ and the elements in $H_i$.
Drutu-Sapir \cite[Theorem 8.6]{DS} proved 
the relative hyperbolicity in the above sense
implies the weak relative hyperbolicity.

\subsection{Toral relatively hyperbolic groups}

We generalize Theorem \ref{krat} to  toral relatively hyperbolic groups.
\begin{thm}\label{main}
Assume $(G, \mathcal H)$ is a  relatively hyperbolic group such that 
 for each $H_i \in \mathcal H$, $H_i$ contains $\Z^{n_i}$ as a finite index subgroup.
Assume $G$ acts on $X$ properly and co-boundedly by isometries
for geodesic metrics $d_1, d_2$ (or more generally, $d_2$
is a coarsely geodesic metric). 
If they are weakly asymptotically isometric, then 
they are coarsely equal.
\end{thm}


Since $d_2$ may be coarsely geodesic, 
the theorem for example applies to two word metrics on $G$
(take the Cayley graph for one of the two metrics as $X$).
We assume $d_1$ is geodesic since we apply \cite{DS} and \cite{sisto}
to $d_1$. They only discuss geodesic metrics. It looks  likely that the results we use
from those two papers hold for coarsely geodesic metrics.

Notice that $\Gamma(G,S)$ and $X$ are quasi-isometric with respect to both $d_1, d_2$
(\cite[Proposition 8.19]{bridson}. The statement therein is
for co-compact group actions on a geodesic space, but the argument uses only co-boundedness, and also it applies to a coarse geodesic space).

In the proof of the theorem we will  first apply
the theorem by Burago \cite{burago} to 
$(H,d_1|H), (H, d_2|H)$, and 
want to conclude $d_1 \sim d_2$ on $H$.
A little issue is that $d_i|H$ is not  a geodesic metric.
So we first modify his result then argue that 
the modification is enough for us.

\begin{thm}\label{burago}
Let $d$ be a coarsely geodesic metric on $\Z^n$
that is invariant by the left action of $\Z^n$.
Then 
\begin{enumerate}
\item 
for each $g \in \Z^n$, 
$\lim_{n \to \infty} \frac{d(1,g^n)}{n}$ exists.
We write the limit by $|g|_d$.
\item
There is a constant $C$ such that for all $g$, 
$|d(1,g) - |g|_d| \le C$.
\item
$|d(1,g^n) - n |g|_d| \le C$ for all $n>0$ and $g$.
\end{enumerate}
\end{thm}

\proof
We first prove that the limit exists. 
The argument is a modification of 
the proof of \cite[Theorem 1]{burago}.
We explain the change we need. 
We only need to modify Lemma 4 in \cite{burago}.

Embed $\Z^n < \R^n$ as a subgroup. Fix $x \in \Z^n$.
The claim of Lemma 4 is for $h \in \Z^n$,
$2d(x, h.x) \sim d(x, h^2(x))$, where a bound does
not depend on $h$. To prove it, 
 join $x, h^2.x$ by a coarse
geodesic $\gamma:[0,L] \to \Z^n$. Approximate it by a 
continuous path $\gamma':[0,L] \to \R^n$ such that at each time,
the points on the two paths are boundedly apart such that 
for each point $g=\gamma(t)$ on $\gamma$, there exists
$t'$ such that $\gamma(t')=\gamma'(t')$ (also 
$|t-t'|$ is bounded). 
In other words, $\gamma'$ and $\gamma$ visit the same 
points in $\Z^n$ at the same time (for each point).

Apply Lemma 2 in \cite{burago} to $\gamma'$ 
and divide it into at most $n$ segments at points in $\Z^n$ and rearrange, then get a path $\gamma'':[0,L]\to \R^n$ 
from $x$ to $h^2.x$  such 
that the distance from $h.x$ to $\gamma''$
is bounded. In the  original case it exactly passes $h.x$ since 
we can divide the path anywhere. In our setting, we 
approximate the original dividing points by nearby points in $\Z^n$.
Since $\Z^n$ is commutative and we cut only at most 
$n$ times, we get a uniform bound. 
From $\gamma''$, by approximating it by nearby
points in $\Z^n$, we get a coarse geodesic $\gamma''':
[0,L] \to \Z^n$ from $x$ to $h^2(x)$.
It follows that 
$2d(x, h.x)=d(x, h.x) + d(h.x, h^2.x) \sim h(x, h^2(x)).$
Lemma 3 in \cite{burago} is proved, therefore (1) and (2) are  proved
as in \cite{burago}.

For (3), notice that by the definition of $|g|_d$, 
we have $|g^n|_d = |n||g|_d$, therefore, 
$|d(1,g^n) - n |g|_d| \le C$ for all $n>0$ and $g$.
\qed

\begin{cor}\label{burago2}
Let $d_1, d_2$ be a coarsely geodesic metric on $\Z^n$
that is invariant by the left action of $\Z^n$.
Assume that if $d_2(1,g^n)$ or $d_1(1,g^n) \to \infty$
as $n \to \infty$, then 
$\lim_{n \to \infty} \frac{d_1(1,g^n)}{d_2(1,g^n)}=1$.
Then $d_1 \sim d_2$.

Moreover, the result is true if we replace $\Z^n$
by a group which contains $\Z^n$ as a finite index 
subgroup.
\end{cor}

In other words, if $d_1, d_2$ are weakly asymptotically isometric,
then they are coarsely equal. 
\proof
By Theorem \ref{burago}, there is $C$ such that 
for all $n>0, g$, we have 
$|d_1(1,g^n) - n |g|_{d_1}| \le C$ and 
$|d_2(1,g^n) - n |g|_{d_2}| \le C$.
Now by our assumption, $|g|_{d_1}=|g|_{d_2}$.
Again, by the theorem, for all $g$,
$|d_1(1,g) - d_2(1,g)|\le 2C$. Since both $d_1, d_2$
are $\Z^2$ invariant, we have $d_1 \sim d_2$.

For the moreover part, we already know
$d_1 \sim d_2$ on $\Z^n$.
But any element in $H$ is at bounded distance, say $D$, from 
a subgroup isomorphic to $\Z^n$
(for both $d_1, d_2$), therefore $d_1\sim_{C+2D} d_2$
on $H$.
\qed


Now, here is a lemma that will assure that $d|_{H_i}$ is coarsely geodesic.
\begin{lemma}\label{coarse}
Assume $H$ acts on a coarsely geodesic space $X$ with a point 
$x\in X$.
Define a metric on $H$ by $d(a,b)=d_X(a.x, b.x)$.
Assume that there is $C$ such that for any $a,b \in H$, 
there is a $X$-geodesic between $a.x, b.x$ which is in the 
$C$-neighborhood of $H.x$. 
Then $d$ is a coarse geodesic metric on $H$.

\end{lemma}
\proof
This is straightforward from the definition.
First assume that $X$ is geodesic.
Let $\gamma(t), 0 \le t \le L$ be a geodesic
in $X$ from $a.x$ to $b.x$.
For each $t$, choose $h_t \in H$ with 
$d(\gamma(t),h_t.x) \le C$.
Define a path $\alpha$ in $H$ from $a$ to $b$
by $\alpha(t)=h_t, 0 \le t \le L$.
Now for any $0 \le t \le  s \le L$, we have 
$d(\alpha(t), \alpha(s))=d(h_t,h_s) =d_X(h_t.x, h_s.x)
\sim_{2C}d_X(\gamma(t), \gamma(s)) =|t-s|$.
Therefore $\alpha$ is a $2C$-coarsely geodesic. 
If $X$ is coarsely geodesic, then start with a 
coarse geodesic $\gamma$ and argue.
\qed

%

We quote several results on asymptotically tree graded spaces and relatively hyperbolic 
groups.
First, a peripheral subgroup is almost convex, therefore
Lemma \ref{coarse} applies to $d|_H$.

\begin{lemma}\label{almost.convex}\cite[Lemma 4.3]{DS}
Let $(G,\mathcal H)$ be relatively hyperbolic. Then each
$H \in \mathcal H$ is almost convex in $G$, in the 
sense that any $(K,L)$-quasi-geodesic joining
two points of $H$ in $\Gamma(G,S)$ is in the $C$-neighborhood of $H$,
where $C$ depends on $K, L$ but not on the quasi-geodesic.

\end{lemma}
We remark that in \cite{DS} the lemma is stated only for geodesics, which is 
sufficient to apply Lemma \ref{coarse}, but the lemma 
holds for quasi-geodesics. The reason is that two distinct
$aH$ and $bH'$ stay close only in a bounded set (Lemma \ref{disjoint}), so that the claim follows using the next lemma. 

Two quasi-geodesics in $X$ with common end points
stay close to each other in the following sense.
We say that the two quasi-geodesics {\it fellow travels}.
This is a version of  Morse lemma for asymptotically 
tree graded spaces.

\begin{lemma}\label{bcp}\cite[Theorem 1.12]{DS}
Let $x,y \in X$, and $\alpha$ a $(K,L)$-quasi-geodesic 
and $\gamma$ a geodesic both between $x,y$ in $(X,d_1)$. Then
there exists a constant
$C(K,L)$ such that they $C$-fellow travel up to $aH.x$'s,
($H \in \mathcal H$).
More precisely, $\alpha$ is in the $C$-neighborhood of $\gamma$
except for the union of some (long) sub-quasi-geodesics of $\alpha$ each of which 
is contained in the $C$-neighborhood of some $aH.x$
such that $aH.x$ is distance at most $C$ from $\gamma$.
In this case, the end points of  each of the sub-quasi-geodesics
are in the $C$-neighborhood of $\gamma$.

\end{lemma}
In \cite{DS} the above lemma is stated for $\Gamma(G,S)$, but $\Gamma(G,S)$
 and $(X,d_1)$, as well as $(X,d_2)$, are quasi-isometric
therefore the results hold for $(X,d_1), (X,d_2)$ as well. 
%
%
%
%

In the above lemma, we may assume the long subpaths of $\alpha$
are disjoint by the following lemma.

\begin{lemma}\label{disjoint}\cite[Lemma 4.7]{DS}
Fix $x \in X$. For each $D$, the diameter of the intersection 
of the $D$-neighborhood of $aH_i.x$ and the $D$-neighborhood
of $bH_j.x$ in $X$ is uniformly bounded  unless $aH_i=bH_j$.
The bound
depends only on $D$.

\end{lemma}

Following \cite{DS}, the (almost) {\em projection} to $aH$,
$\pi_{aH}$, in $\Gamma(G,S)$ is defined  as follows for $H \in \mathcal H$:
for $g \in G$, $\pi_{aH}(g)$ is the subset of points in $aH$
whose distance from $g$ is less than $d(g,aH)+1$.

The following result also holds for the projection to $aH.x$ in $X$
as well (the proof is same).

\begin{lemma}\label{projection}\cite[Lemma 1.13 (1) and Theorem 2.14]{sisto}
Let $(G,\mathcal H)$ be relatively hyperbolic.
Let $H \in \mathcal H$ and $\pi$ be the  projection 
to $aH$ 
in the Cayley graph with $a \in G$. Then
\begin{enumerate}
\item
Any $(K,L)$-quasi geodesic from a point $g$ in $\Gamma(G,S)$
to a point in $aH$ passes the $C$-neighborhood of 
$\pi(g)$. The constant $C$ depends on $K, L$, but 
not on $a, H$ and $g$.
\item
The diameter of $\pi(g)$ is uniformly bounded. The bound does not
depend on $a, H$ and $g$.
\end{enumerate}
\end{lemma}

We may write $\Gamma(G, S \cup \mathcal H)$
as $G'$. We denote the distance on $\Gamma(G, S \cup \mathcal H)$ by $d_{G'}$.
$\Gamma(G, S \cup \mathcal H)$ is hyperbolic (see \cite[Section 8]{DS})
 and $G$ acts on it. 
Each edge in $\Gamma(G, S \cup \mathcal H)$ that is not in $\Gamma(G,S)$ 
joins two points 
in $aH$ for some $H \in \mathcal H$.
Given a geodesic $\gamma$ in $\Gamma(G, S \cup \mathcal H)$, a {\it lift} is a path in $\Gamma(G,S)$
obtained by replacing each edge of $\gamma$ that is not in $\Gamma(G,S)$
by a geodesic in $\Gamma(G,S)$ connecting the two end points of the edge.

\begin{lemma}\label{lift.geodesic}\cite[Prop 1.14]{sisto}
If $\gamma$ is a geodesic in $\Gamma(G, S \cup \mathcal H)$ then its
lift is a quasi-geodesic in $\Gamma(G,S)$ with uniform quasi-geodesic 
constants.

\end{lemma}
The above lemma does not hold for quasi-geodesics
$\gamma$ in general. 

\begin{lemma}\label{projection.geodesic}\cite[Prop 8.25]{DS}
Let $\gamma$ be a quasi-geodesic in $\Gamma(G,S)$
between $x,y$ and $\gamma'$ a quasi-geodesic
in $\Gamma(G, S \cup \mathcal H)$ between $x,y$.
Then they are in a bounded Hausdorff-distance in $\Gamma(G, S \cup \mathcal H)$.
The bound depends only on the quasi-geodesic constants.

\end{lemma}

%
%

We start the proof of Theorem \ref{main}.
\proof
Fix $x \in X$. Define a left invariant (pseudo-)metric $d_1(g,h)=d_1(g.x, h.x)$ 
on $G$, and also $d_2$ in the same way.
Both $(G,d_1), (G,d_2)$ are quasi-isometric to $\Gamma(G,S)$.
As before set $\D=d_1-d_2$ on $X$ and $G$. Note that 
$\D(g,h)=\D(g.x,h.x)$.

We summarize what we know on each peripheral subgroup $H$ by now. 
On $H$, $d_1$ and $d_2$ are coarsely geodesic metrics by Lemma \ref{coarse}
and Lemma \ref{almost.convex}.
By assumption they are weakly asymptotic on $H$.
Therefore $d_1 \sim d_2$ on $H$ by Corollary \ref{burago2}.
The conclusion holds for $aH$ as well. 

%


\begin{lemma}\label{constant}
There exists $L$ such that on each $aH.x$ and $z \in X$, 
$\D(z,y)$ varies at most $L$ for $y \in aH.x$. $L$ does
not depend on $a, H$ and $z$.

Moreover, the statement holds when $y$ is in the $K$-neighborhood
of $aH.x$  ($L$ depends  on $K$). 

Also, the statement holds on $G$, namely,
on each $aH$ and $g \in G$, $\D(g,h)$ varies at most $L$ for $h$
that is in the $K$-neighborhood of $aH$. 

\end{lemma}

\proof
As usual we assume $d_2$ is geodesic. We omit details
for the coarse geodesic case. 
Fix a point $A \in \pi_{aH.x}(x)$, where $\pi_{aH.x}$
is defined in $(X,d_1)$.
We claim that $\D(x,y) \sim \D(x,A)$
for all $y \in aH.x$ such that the constant for $\sim$ does
not depend on $y$ or $aH$.
Let $\gamma_1, \gamma_2$ be geodesics from $x$ to $y$
for $d_1, d_2$. Notice that $\gamma_2$
is a $d_1$-quasi-geodesic with controlled constants.
By Lemma \ref{projection}, they pass $C$-neighborhood (in both $d_1$ and $d_2$,
which are quasi-isometric to each other)
 of $A$.
$C$ does not depend on $a, H, y$.
Now take $q_1, q_2$ on each geodesic with $d_i(q_1,A) \le C$ and $d_i(q_2,A) \le C$ with $i=1,2$.
Then, by Lemma \ref{basic.additive} and Lemma \ref{triangle}, we have  
$\D(x,y) \sim_{4C} \D(x,q_1) + \D(q_1, y)
\sim_{8C} \D(x,A) +\D(A,y)$.

But we already know $d_1\sim d_2$ on $aH$, therefore, we have $\D(A,y) \sim 0$.
It follows $\D(x,y) \sim \D(x,A)$.
The argument is complete when $y \in aH.x$.
The moreover part now follows from Lemma \ref{triangle}.
The argument for the group $G$ is same. 
\qed

Here is a consequence. 
\begin{lemma}\label{estimate}
There exists $P$ such that for 
any  $x,y,z \in G$, we have 
$$|\D(x,y)-\D(x,z)| \le P d_{G'}(y,z).$$

\end{lemma}
\proof
Let $N=d_{G'}(y,z)$ and $\gamma$ a geodesic from $y$ to $z$ in $\Gamma(G, S \cup \mathcal H)$ with 
vertices $y=y_0, y_1, \cdots, y_N=z$.
$y_n$ and $y_{n+1}$ are joined by an edge.
There exists $K$ such that $|\D(x,y_n)-\D(x,y_{n+1})| \le K$ 
if the edge is in $\Gamma(G,S)$.
On the other hand, if the edge is not in $\Gamma(G,S)$, then 
$y_n, y_{n+1}$ are in some $aH$, therefore $|\D(x,y_n)-\D(x,y_{n+1})| \le L$ 
by Lemma \ref{constant}.
Now it follows that $|\D(x,y)-\D(x,z)| \le (L+K)N$.
Set $P=L+K$.
\qed

\begin{lemma}\label{subadditive}
Let $\gamma$ be a $d_2$-quasi-geodesic in $X$.
 Then, $\D$ is almost additive on $\gamma$.
Namely, let $x,z,y$ be points on $\gamma$ in this order, then 
$\D(x,z)+\D(z,y) \sim_B \D(x,y)$,
where $B$ depends on the quasi-geodesic constants of $\gamma$.

The statement holds for quasi-geodesics 
on $\Gamma(G,S)$ as well. 
\end{lemma}

\proof
By assumption $\gamma$ is a $d_2$-quasi-geodesic with controlled
constants. Then it is a $d_1$-quasi-geodesic with controlled
constants as well. 
Let $\gamma'$ be a $d_1$-geodesic in $X$ from $x$ to $y$.
By Lemma \ref{bcp}, there exists $C$ that depends on the 
quasi-geodesic constants such that 
$\gamma$ and $\gamma'$ stay close to each other 
except for subsegments in $\gamma$ each of which stays
in the $C$-neighborhood of one $aH.x$ for a long time, but 
the end points of those segments are $C$-close to $\gamma'$.
If $z$ is outside of those segments, then 
$z$ is close to $\gamma'$, which implies
the almost additivity at $z$ by Lemma \ref{basic.additive}.


Now assume $z$ is contained in one of the subsegments, say,
$[z_1,z_2]$.
Each $z_i$ is $C$-close to $\gamma'$, therefore  by Lemma \ref{basic.additive},
we have
$\D(x,y) \sim \D(x,z_1) + \D(z_1,y)$. 
On the other hand, since $z,z_1,z_2$ is in the $C$-neighborhood of $aH.x$,
by Lemma \ref{constant} and Lemma \ref{triangle}, 
we have 
$ \D(x,z_1) + \D(z_1,y) \sim \D(x,z) + \D(z,y)$
(Lemma  \ref{constant} applies  to  $\D(*,y)$). Combining them, $\D$ is almost additive at $z$.
The argument is complete. 

Since $\Gamma(G,S)$ and $X$ are quasi-isometric, we also have
the almost additivity in $\Gamma(G,S)$.
\qed

We go back to the proof of the theorem. 
We want to show that $\D(x,y)$ is bounded on $X$, which 
 is equivalent to that $\D(g,h)$ is bounded on $G$.
To argue by contradiction, assume not. We will find $f \in G$
such that $\D(x, f^n.x)=\D(1,f^n)$ grows roughly linearly on $n$, which 
will be a contradiction since $d_1$ and $d_2$ are weakly asymptotic.

Take $g$ such that $\D(1, g)$ is very large.
Let $\gamma$ be a geodesic from $1$ to $g$ in $\Gamma(G,S)$, and 
let $h \in G$ be such that $h$ is on $\gamma$ and 
$\D(1,h)$ is approximately $\D(1,g)/2$.
This is possible since $\D(1,y)$ is almost continuous 
when we vary $y$ on $\gamma$.

Let $\alpha$ be a geodesic from $1$ to $g$ in $\Gamma(G, S \cup \mathcal H)$.
By Lemma \ref{projection.geodesic}, the Hausdorff distance 
between $\gamma$ and  $\alpha$ in $\Gamma(G, S \cup \mathcal H)$ is bounded. 
In particular, $h$ is at bounded distance from  $\alpha$ in $\Gamma(G, S \cup \mathcal H)$.
Notice that $d_{G'}(1,h), d_{G'}(h, g), d_{G'}(1, g)$ are all large
since $\D(1,h), \D(h,g), \D(1,g)$ are all large (use Lemma \ref{estimate}).
Since $\Gamma(G, S \cup \mathcal H)$ is hyperbolic, as before, 
one of the paths $h^*, k^*, (hk)^*$ (we use the same 
notation as in the proof of Theorem \ref{krat} with $x=1$) is a quasi-geodesic in $\Gamma(G, S \cup \mathcal H)$
with uniform quasi-geodesic constants.
Denote it by $\gamma=f^*$.

First, we assume that $\gamma$ is a geodesic, and argue.
Take a lift, $\beta$, of $\gamma$ in $\Gamma(G,S)$. 
Then $\beta$ is a quasi-geodesic in $\Gamma(G,S)$ with controlled
constants  by Lemma
\ref{lift.geodesic}. 
Since $\D(1,f)$ is very large (approximately 
$\D(1, g)$ or $\D(1, g)/2$), by Lemma \ref{subadditive}, 
$\D(1, f^n)$ grows roughly linearly on $n$
(apply the lemma at each point $f^n$).

In general, $\gamma$ is only a quasi-geodesic in $\Gamma(G, S \cup \mathcal H)$ with the quasi-geodesic
constants controlled. In this case,  for each $N>0$, take a
geodesic $\gamma'$ from $1$ to $f^N$ in $\Gamma(G, S \cup \mathcal H)$,
 so that $1,f, \cdots, f^{N}$ are in a 
bounded neighborhood of $\gamma'$ in $\Gamma(G, S \cup \mathcal H)$. This is 
because $\Gamma(G, S \cup \mathcal H)$
is hyperbolic. For each $f^n$, let $y_n \in \gamma'$ be a closest
point  on  $\gamma'$ for $d_{G'}$. Take $y_0=1, y_N=f^N$.
By Lemma \ref{estimate}, for each $n$, $\D(y_n, y_{n+1}) \sim \D(f^n, f^{n+1})
= \D(1, f)$.
Take a lift of $\gamma'$, denoted by $\beta$, which is a quasi-geodesic in
 $\Gamma(G,S)$
with controlled constants. 
The points $y_n$ are on $\beta$, and by Lemma \ref{subadditive},
$\D(y_0, y_n)$ grows roughly linearly on $n$
(roughly the slope is $\D(1, f)$, which is much larger than any constants
for $\sim$ in the above argument).
Again by Lemma \ref{estimate}, $\Delta(1,f^n)$ grows roughly linearly 
for $0 \le n \le N$, with the slope roughly 
$\D(1, f)$.
Since $N$ was arbitrary, $\D(1,f^n)$
grows roughly linearly on $0 \le n$.
This finish the argument. 
\qed

\subsection{More general case}\label{complex}
In the proof of Theorem \ref{main}, what we need from 
the peripheral subgroups $H$ is the property that any two 
weakly asymptotic, coarsely geodesic metrics on $H$ are
coarsely equal (see the discussion in the beginning 
of the proof). We verified this property in Corollary \ref{burago2}
for virtually abelian groups. 
We restate Theorem \ref{main} in this more general form.
The proof is identical and we omit it.

 \begin{thm}\label{main2}
Let $(G, \mathcal H)$ be a  relatively hyperbolic group.
Assume that each $H_i \in \mathcal H$ satisfies the property 
such that  any two 
weakly asymptotic, coarsely geodesic metrics on $H_i$ are
coarsely equal. 
Suppose  $G$ acts on $X$ properly and co-boundedly by isometries
for geodesic metrics $d_1, d_2$ (or more generally, $d_2$
is a coarsely geodesic metric). 
If $d_1$ and $d_2$ are weakly asymptotically isometric on $X$, then 
they are coarsely equal.
\end{thm}

One potential application would be to non-uniform lattices
in the Lie group $SU(n,1)$. 
 It is known that such lattice
is relatively hyperbolic with peripheral subgroups
virtually nilpotent. As we mentioned in the introduction,
the assumption in the theorem does not hold
in general for all nilpotent groups.

\section{Marked length spectrum}\label{marked}
In this section we discuss an application to the marked length spectrum 
problem. 
We state a variant of Theorem \ref{main}.

\begin{thm}\label{rel.hyp}
Let $(X_1, d_1), (X_2,d_2)$ be 
geodesic metric spaces  on which $G$ acts
by isometries, co-boundedly and properly with respect to both $d_1$ and $d_2$. 
Suppose the action on $X_1$ is free. 
We allow $X_1$ to be coarsely geodesic.  
Assume 
$|g|_1=|g|_2$ for any hyperbolic element $g \in G$.
If $G$ is toral relatively hyperbolic group, then there exists a $G$-equivariant,
$(1,C)$-quasi-isometry map 
$f:X_1\to X_2$. 

\end{thm}
\proof
We follow the argument for Theorem \ref{main}.
We fix points $x_1 \in X_1$ and $x_2 \in X_2$, 
and define $\D(g,h)=d_1(g.x_1, h.x_1)-d_2(g.x_2, h.x_2)$
for $g,h \in G$.
We claim that  $\D$ is bounded. We will prove this later,
but once this is known, then 
there exists a desired map $f$. In other words,
$\D(p,q)=d_1(p,q)-d_2(f(p),f(q))$ is bounded.
Indeed, set $f(g.x_1)=g.x_2$
for each $g \in G$ (use the $G$-action on $X_1$ is free). 
Then, $f$ is $G$-equivariant and $\D(p,q)$ is bounded for $p,q \in G.x_1$.
Moreover, we can extend $f$ to $X_1$, $G$-equivariantly,  
such that $\D$ is bounded. 
Indeed, choose a point $p$ on each $G$-orbit in $X_1$. There is $g \in G$ such that 
$d_1(p,g.x_1)$ is bounded. Define $f(p)=g.x_2$, and extend $f$,
$G$-equivariantly, to the $G$-orbit of $p$.
It is clear that $\D$ is bounded on $X_1$.

We are left to show $\D$ is bounded on $G.x_1$.
Again, we repeat the argument for  Theorem \ref{main}.
When needed, we replace $X$ with $X_2$ and map any objects
in $X_1$ to $X_2$ by $f$ and argue on $X_2$.
Notice that the map $f$ defined above is a quasi-isometry from $X_1$
to $X_2$.
For example, in the proof of  Lemma \ref{constant}, 
$\gamma_1$ is a geodesic from $x_1$ to $y_1=ah.x_1 \in aH.x_1$, and 
$\gamma_2$ is a geodesic from $x_2$ to $y_2 =ah.x_2\in aH.x_2$.
Then $f(\gamma_1)$ is a quasi-geodesic from $x_2$ to $y_2$.
Apply the argument to $\gamma_2$ and $f(\gamma_1)$.
Lemma \ref{subadditive} is similar. 
We omit details. 
\qed

We apply the result to the marked length spectrum problem for manifolds.

\begin{cor}\label{manifold}
Let $(M_1,d_1), (M_2,d_2)$ be closed Riemannian manifolds
with the isomorphic fundamental group $G$ that 
is toral relatively hyperbolic.
Assume they have the same marked length spectrum.
Then there is a $G$-equivariant $(1,C)$- quasi-isometry 
$f:\tilde M_1 \to \tilde M_2$.

Moreover, if there is a homeomorphism $H: M_1 \to M_2$
that induces the isomorphism on $G$ and  lifts
to a  $G$-equivariant homeomorphism $h:\tilde M_1 \to \tilde M_2$,
then $h$ is a $(1,C')$- quasi-isometry.

\end{cor}


\proof
Let $X_i$ be the universal cover of $M_i$, respectively. 
Each of them have the lift of $d_i$, 
which we also write by $d_i$. 
Each action by $G$ on $X_i$ is free. By assumption, we 
can apply Theorem \ref{rel.hyp} and 
we get a desired map $f$.

For the moreover part, it suffices to show that there
is a constant $L$ such that for any point $p \in \tilde M_1$,
$d_2(f(p),h(p)) \le L$.
To see this, fix $x_1 \in X_1$ and set $x_2=h(x_1)$.
Let $\pi:X_1 \to M_1$ be the covering map.
Given $p \in X_1$, join $\pi(x_1)$ to $\pi(p)$ by a shortest
geodesic $\gamma$ in $M_1$.
Lift it to a geodesic $\tilde \gamma$ from $p$ to $g.x_1$
with $g \in G$.
Then $h(\tilde \gamma)$
is a path from $h(p)$ to $h(g.x_1)=g.(h(x_1))=g.x_2$.
By the argument for Theorem \ref{rel.hyp}, $f(p)=g.x_2$ (to be precise,
we can extend $f$ 
in this way).
That means that $h(p)$ and $f(p)$ are joined by $h(\tilde \gamma)$, 
but the length of this path is bounded since the length of $\tilde \gamma$
is bounded and $h$ is continuous and $G$-equivariant. 
\qed

%


For example, if $M_1$ has non-positive curvature,
dimension is not 3 nor 4, and $M_2$ is aspherical
(for example $M_2$ has non-positive curvature), 
then each isomorphism of $G$ is induced 
by a homeomorphism $H$ by Farrell-Jones (see \cite{farrell}), 
therefore its lift $h$ is a $(1,C)$-quasi-isometry map. 

The case where $G$ is word-hyperbolic in Corollary \ref{manifold} is proved by Furman \cite[Theorem 2]{furman}.
His argument is different from ours, and 
uses Patterson-Sullivan measures for hyperbolic 
groups with respect to word metrics constructed by Coornaert \cite{co}, 
and does not seem to apply to prove Theorem \ref{rel.hyp}.

If $C=0$ for some $f$ in Theorem \ref{manifold}, $M_1$ and $M_2$ are isometric.
This is the conclusion that {\it the marked length spectrum problem/conjecture}
concerns (see the conjecture in \cite[3.1]{katok}
for negative curvature case). Several cases are known to have 
the positive answer, 
for example for surfaces of negative curvature \cite{otal}.
Interestingly, Bonahon \cite{bonahon} 
gave examples of 
geodesic metrics $d_1, d_2$ on a hyperbolic surface
such that $d_1$ is Riemannian with constant negative curvature 
and that they have the same marked length spectrum, but are 
not isometric to each other. Theorem \ref{rel.hyp}
applies to his example. 
For higher dimension there is a result by Hamenst\"adt
($M_1$ is a rank-1 locally 
symmetric space and $M_2$ is negatively curved) 
using a theorem of Besson-Courtois-Gallot \cite{bcg} on the volume entropy.
In all of those cases, the fundamental 
group is hyperbolic. 
Our result  put  the marked length spectrum 
problem (for manifolds) into context for a broader class
of groups.
Although this is not for manifolds, 
another case where 
the isospectral length problem is solved is 
$\R$-trees under the assumption that the action 
is minimal and  semi-simple, \cite{culler}.

\end{document}